\documentstyle[art11]{article}
\font\Bbb=msbm10 at 11pt
\font\bbb=msbm10 at 8pt
\font\BBb=msbm10 at 14pt
\font\aaa=msam10 at 11pt
\def\ppp{\mbox{\Bbb Y}}
\def\real{\mbox{\Bbb R}}
\def\integer{\mbox{\Bbb Z}}
\def\mod #1{\,(\mathop{\rm mod}#1)}
\def\wall{{\cal W}}
\def\dbrack#1{[\![#1]\!]}
\newdimen\mypt
\newdimen\halfpt
\def\pic#1 #2 #3;{\raisebox{0.6ex}{\raisebox{-#2\halfpt}{\vbox
      to #2\mypt{\hbox to #3\mypt{\special{em:graph
      #1.bmp}\hfill}\vfill}}}}

\def\proof{\addvspace{\medskipamount}\noindent
      {\sc Proof.} }
\def\proved{\ifmmode\eqno\Box\medskip\else
        \nobreak\hfill\nopagebreak\discretionary{}
        {\hbox to\textwidth{\hfill$\Box$}}{\hbox{$\Box$}}\par
        \addvspace\medskipamount\fi}
\def\alf{\mbox{\Bbb A}}
\def\brack#1{\Big\langle#1\Big\rangle}
\newcounter{defin}
\long\def\definition#1{\medbreak\noindent\refstepcounter{defin}%
{\bf Definition \arabic{defin}. }#1\par\medskip}
\newcounter{remark}
\long\def\remark#1{\medbreak\noindent\refstepcounter{remark}%
{\sl Remark \arabic{remark}. }#1\par\medskip}
\font\msamfive=msam5 at 3pt
\def\uparc(#1_#2#3_#4){{\mathop{#1_{#2}#3_{#4}}\limits^{\displaystyle
\frown\kern-3pt\llap{\raisebox{1.5pt}{\msamfive\char121}}
\phantom{\scriptscriptstyle#4}
}}}
\def\loarc(#1_#2#3_#4){{\mathop{#1_{#2}#3_{#4}}\limits_{\displaystyle
\smile\kern-3pt\llap{\raisebox{1.7pt}{\msamfive\char113}}
\phantom{\scriptscriptstyle#4}
}}}
\def\diez{\mbox{\sc\char35}}
\def\Sing{{\rm Sing \,}}
\def\sing{{\rm sing \,}}

\author{S. D. Tyurina
}
\date{Dep. of Phys. and Math.\\
Kolomensky Pedagogical Institute\\
Kolomna 140411 RUSSIA\\
e-mail: tyurina@mccme.ru}
\title{{\sc Explicit formulas for the Vassiliev knot invariants.}\\
{\scriptsize(The brief version of this paper is in materials of
International Conference "Monodromie et
\'equations diff\'erentielles en th\'eorie des singularit\'es et des
repr\'esentations de groupes." Th\`eses. -- Luminy, France, CNRS-SMF,
jan.-1999.
)}}

\begin{document}
\baselineskip=11.8pt
\maketitle
%
\section{Knots and singular knots.}
Let $K:S^1 \longrightarrow {\bf R^3}$
be an oriented knot and $K^{\sing}_n:S^1
\longrightarrow {\bf R^3}$ be a singular knot
with $n$ double points.

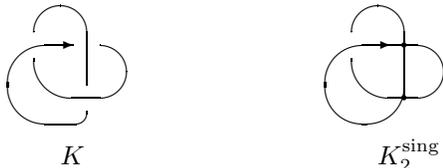
\begin{picture}(400,65)
\put(100,45){\vector(1,0){5}}
\put(100,50){\oval(20,20)[t]}
\put(110,30){\line(0,1){20}}
\put(105,30){\oval(50,30)[tl]}
\put(110,40){\oval(40,30)[bl]}
\put(110,25){\line(1,0){5}}
\put(95,20){\oval(30,10)[br]}
\put(115,35){\oval(20,20)[r]}
\put(95,30){\oval(30,30)[bl]}

\put(220,45){\vector(1,0){5}}
\put(220,45){\line(1,0){15}}
\put(220,50){\oval(20,20)[t]}
\put(230,30){\line(0,1){20}}
\put(225,30){\oval(50,30)[tl]}
\put(230,40){\oval(40,30)[bl]}
\put(230,25){\line(0,1){15}}
\put(230,25){\line(1,0){5}}
\put(215,30){\oval(30,30)[b]}
\put(230,25){\circle*{2}}
\put(230,45){\circle*{2}}
\put(235,35){\oval(20,20)[r]}

\put(100,0){$K$}
\put(220,0){$K^{\sing}_2$}
\end{picture}

Denote by ${\cal K}$ the space of knots and
by ${\cal K}^{\sing}$ the space of singular knots.

\section{The Vassiliev invariants.}
Let $V:{\cal K} \longrightarrow {\bf Q}$ be a knot invariant.
It may be extended from ordinary knots to singular knots. Define
$i$-th derivative of invariant $V \quad $
$V^{(i)}:{\cal K}^{\sing} \longrightarrow {\bf Q}$ inductively
via the Vassiliev skein-relation:
$$V^{(0)} = V,$$

\begin{picture}(420,22)
\put(80,9){$V^{(i)}($}
\put(120,11){\circle*{2}}
\put(132,9){$)=V^{(i-1)}($}
\put(217,9){$)-V^{(i-1)}($}
\put(299,9){).}
\put(109,0){\vector(1,1){22}}
\put(131,0){\vector(-1,1){22}}
\put(193,0){\vector(1,1){22}}
\put(215,0){\line(-1,1){10}}
\put(203,12){\vector(-1,1){10}}
\put(274,0){\line(1,1){10}}
\put(288,12){\vector(1,1){10}}
\put(297,0){\vector(-1,1){22}}
\end{picture}

{\bf Definition.} A knot invariant $V:{\cal K} \longrightarrow
{\bf Q}$ is said to be Vassiliev invariant of order
less than or equal to $n$,
if it's $(n+1)$-st derivative vanishes identically:
$$V^{(n+1)} \equiv 0.$$

We will denote Vassiliev invariant of order less than or equal to $n$ by $V_n$ and the space
of Vassiliev invariants over {\bf Q} by ${\cal V}_n.$
The spaces of Vassiliev invariants have following dimensions:

\begin{center}
\begin{tabular}{|c||c|c|c|c|c|}
\hline
$n$ & $0$ & $1$ & $2$ & $3$ & $4$\\
\hline
$\dim {\cal V}_n$ & $1$ & $1$ & $2$ & $3$ & $6$\\
\hline
\end{tabular}
\end{center}

\section{Chord diagrams.}
A chord diagram $D_n$ of order $n$ is a circle with
a distinguished set of $n$ unordered pairs of points connected by chords
regarded up to orientation preserving diffeomorphisms of the circle.
Denote by ${\cal D}_n$
the space generated by chord diagrams of order $n$ over ${\bf Q}.$ For example,
there exist only two non-equivalent chord diagrams of order 2:
\begin{picture}(85,10)
\put(0,0){${\cal D}_2=\{$}
\put(43,3){\circle{15}}
\put(68,3){\circle{15}}
\put(76,0){$\}$}
\put(40,-3){\line(0,1){12}}
\put(46,-3){\line(0,1){12}}
\put(63,-2){\line(1,1){10}}
\put(73,-2){\line(-1,1){10}}
\put(55,0){,}
\end{picture}
and space ${\cal D}_3$ contains 5 diagrams:
\begin{picture}(160,10)
\put(96,3){\circle{15}}
\put(71,3){\circle{15}}
\put(83,0){,}
\put(46,3){\circle{15}}
\put(58,0){,}
\put(146,3){\circle{15}}
\put(133,0){,}
\put(0,0){${\cal D}_3=\{$}
\put(121,3){\circle{15}}
\put(159,0){$\}.$}
\put(93,-3){\line(0,1){12}}
\put(99,-3){\line(0,1){12}}
\put(116,-2){\line(1,1){10}}
\put(126,-2){\line(-1,1){10}}
\put(108,0){,}
\put(89,3){\line(1,0){15}}
\put(121,-5){\line(0,1){15}}
\put(46,-5){\line(0,1){15}}
\put(43,-3){\line(0,1){12}}
\put(49,-3){\line(0,1){12}}
\put(69,-4){\line(-1,2){5}}
\put(74,-4){\line(1,2){5}}
\put(65,7){\line(1,0){12}}
\put(143,-3){\line(0,1){12}}
\put(148,-4){\line(1,2){5}}
\put(153,-1){\line(-1,2){6}}
\end{picture}

For any singular knot $K^\sing$ we construct it's chord diagram.
A chord diagram of the singular knot is the circle with pre-images
of double points connected with  chords.

\vskip0.2cm
\section{Weight systems.}
{\bf Definition.} A linear function $W:{\cal D}_n \longrightarrow
{\bf Q}$ is called a weight system of order $n$ if it satisfies next relations:

\begin{picture}(420,20)
\put(0,9){1-term relation:}
\put(100,9){$W_n($}
\put(135,11){\circle{15}}
\put(147,9){$)=0$}
\put(135,3){\line(0,1){15}}
\end{picture}

\noindent (this diagram has $n$ chords, one of which is isolated) and

\begin{picture}(420,20)
\put(0,9){4-term relation:}
\put(100,9){$W_n($}
\put(135,11){\circle{15}}
\put(148,9){$)-W_n($}
\put(200,11){\circle{15}}
\put(213,9){$)+W_n($}

\put(265,11){\circle{15}}
\put(280,9){$)-W_n($}

\put(332,11){\circle{15}}
\put(347,9){$)=0$}

\put(134,4){\line(-1,2){5}}
\put(136,4){\line(1,2){5}}

\put(197,4){\line(1,1){10}}
\put(203,4){\line(-1,1){10}}
\put(263,4){\line(1,2){6.5}}
\put(271,14){\line(-1,0){12}}
\put(329,4){\line(1,1){10}}
\put(338,16){\line(-1,0){11.5}}
\end{picture}

\noindent these diagrams have $n$ chords, $(n-2)$ of  which are not drown
here and 2 chords are positioned as shown.

Denote by ${\cal W}_n$ the space of weight systems of order $n$ over {\bf Q}.
The spaces ${\cal W}_n$ have following dimensions:

\begin{center}
\begin{tabular}{|c||c|c|c|c|c|}
\hline
$n$ & $0$ & $1$ & $2$ & $3$ & $4$\\
\hline
$\dim {\cal W}_n$ & $1$ & $0$ & $1$ & $1$ & $3$\\
\hline
\end{tabular}
\end{center}

Examples of weight systems.

\noindent
\begin{picture}(400,40)
\put(0,22){$W_2($}
\put(0,2){$W_2($}
\put(30,25){\circle{15}}
\put(30,7){\circle{15}}
\put(40,22){)=1}
\put(40,2){)=0}
\put(25,20){\line(1,1){10}}
\put(35,20){\line(-1,1){10}}
\put(25,1){\line(0,1){11}}
\put(35,1){\line(0,1){11}}
\put(200,22){$W_3($}
\put(200,2){$W_3($}
\put(230,25){\circle{15}}
\put(230,7){\circle{15}}
\put(240,22){)=2}
\put(240,2){)=1}
\put(225,20){\line(1,1){10}}
\put(235,20){\line(-1,1){10}}
\put(225,1){\line(0,1){11}}
\put(235,1){\line(0,1){11}}
\put(222,25){\line(1,0){15}}
\put(222,7){\line(1,0){15}}
\end{picture}

\hskip6cm $W_3=0$ in other cases.

\noindent
\vskip0.5cm
\begin{tabular}{|c||c|c|c|c|c|c|c|}
\hline
 &
\begin{picture}(20,15)
\put(8,8){\circle{15}}      
\put(6,1){\line(-1,2){5}}   
\put(1,4){\line(1,2){5.5}}     
\put(10,1){\line(1,2){5}}    
\put(15,4){\line(-1,2){5.5}}   
\end{picture}
&
\begin{picture}(20,15)
\put(8,8){\circle{15}}      
\put(5,1.2){\line(0,1){13}}   
\put(11,1.2){\line(0,1){13}}   
\put(1,8){\line(1,0){15}}   
\put(8,0){\line(0,1){15}}   
\end{picture}
&
\begin{picture}(20,20)
\put(8,8){\circle{15}}     
\put(5,1.2){\line(0,1){12.5}}  
\put(11,1.2){\line(0,1){12.5}}  
\put(1,5){\line(1,0){14}}  
\put(1,11){\line(1,0){14}} 
\end{picture}
&
\begin{picture}(20,15)
\put(8,8){\circle{15}}     
\put(8,0){\line(-1,2){6}}  
\put(2,3.4){\line(1,2){6}}   
\put(4,2){\line(3,1){11}}  
\put(4,15){\line(2,-1){12}}
\end{picture}
&
\begin{picture}(20,15)
\put(8,8){\circle{15}}     
\put(8,0){\line(-1,2){6}}  
\put(2,4){\line(1,2){5.5}}   
\put(11,1.52){\line(0,1){12.5}}  
\put(0.51,8){\line(1,0){15}}  
\end{picture}
&
\begin{picture}(20,15)
\put(8,8){\circle{15}}     
\put(5,1.2){\line(0,1){12.5}}  
\put(11,1.2){\line(0,1){12.5}}  
\put(1,5){\line(2,1){14}}  
\put(1,11){\line(2,-1){14}}
\end{picture}
&
\begin{picture}(20,15)
\put(8,8){\circle{15}}     
\put(3,3){\line(1,1){10}}  
\put(13,3){\line(-1,1){10}} 
\put(8,0){\line(0,1){15}}  
\put(1,8){\line(1,0){15}}  
\end{picture}
\\
\hline
$W^1_4$ & $0$ & $0$ & $1$ & $0$ & $0$ & $1$ & $1$\\
\hline
$W^2_4$ & $1$ & $-1$ & $2$ & $0$ & $-1$ & $1$ & $0$\\
\hline
$W^3_4$ & $0$ & $1$ & $-3$ & $1$ & $2$ & $-2$ & $0$\\
\hline
\end{tabular}

\vskip0.2cm
\noindent
Elsewhere $W_4^1=W_4^2=W_4^3=0.$
\section{Symbol of Vassiliev invariant.}
{\bf Definition.} The symbol of Vassiliev invariant of order
$n$ is the restriction of this invariant to singular knots with exactly $n$
double points:
$$\nabla^nV_n = V_n|_{{\cal K}^\sing_n}.$$
The symbol of Vassiliev invariant of order $n$ \quad $\nabla^nV_n$ depends
only on chords diagrams of singular knots. So the symbol is a linear
function on the space ${\cal D}_n:$
$$\nabla^nV_n:{\cal D}_n \longrightarrow {\bf Q}.$$
One can prove that the symbol of Vassiliev invariant $\nabla^nV_n$
satisfies 1-term and 4-term relations.

{\bf Kontsevich theorem.}
{\it Over {\bf Q} each weight system $W_n$ is the symbol of an appropriate Vassiliev
invariant $V_n$
and there exist next exact sequence}
$$0 \longrightarrow {\cal V}_{n-1} \longrightarrow
{\cal V}_n \longrightarrow {\cal W}_n \longrightarrow 0,$$
$${\cal V}_n / {\cal V}_{n-1} \simeq {\cal W}_n.$$

From the Kontsevich theorem follows that we have

-- one basic Vassiliev invariant of order 2,

-- one basic Vassiliev invariant of order 3,

-- three basic Vassiliev invariants of order 4.

\noindent
And it is easy fact that Vassiliev knot invariants of order 0 and 1 are the
constant maps.

\section{Diagrams of knots.}
Let $K$ be a knot with marked base-point $a$ and $x$ be double point of it's
planar projection.

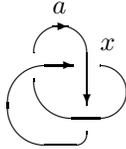
\begin{picture}(400,70)
\put(100,45){\vector(1,0){5}}
\put(100,60){\circle*{2}}
\put(97,60){\vector(1,0){6}}
\put(97,65){$a$}
\put(115,51){$x$}
\put(100,50){\oval(20,20)[t]}
\put(110,50){\vector(0,-1){20}}
\put(105,30){\oval(50,30)[tl]}
\put(110,40){\oval(40,30)[bl]}
\put(110,25){\line(1,0){5}}
\put(95,20){\oval(30,10)[br]}
\put(115,35){\oval(20,20)[r]}
\put(95,30){\oval(30,30)[bl]}
\end{picture}

Numerate branches in neighbourhood of $x$ according to the order of their
passing. Define function $\delta_x$ by next rule:

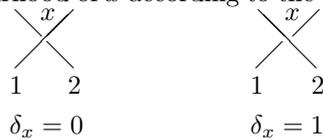
\begin{picture}(400,32)
\put(143,10){\line(1,1){22}}
\put(165,10){\line(-1,1){10}}
\put(141,0){$1$}
\put(163,0){$2$}
\put(141,-15){$\delta_x = 0$}
\put(153,22){\line(-1,1){10}}
\put(152.5,28){$x$}
\put(234,10){\line(1,1){10}}
\put(248,22){\line(1,1){10}}
\put(245,28){$x$}
\put(257,10){\line(-1,1){22}}
\put(232,-15){$\delta_x = 1$}
\put(255,0){$2$}
\put(232,0){$1$}

\end{picture}

\vskip0.5cm
(the orientations of branches are not important).

Define function $\varepsilon_x$ as follows:

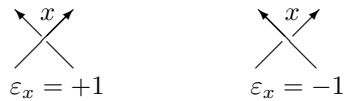
\begin{picture}(400,32)
\put(143,10){\vector(1,1){22}} 
\put(165,10){\line(-1,1){10}}
\put(153,22){\vector(-1,1){10}}
\put(234,10){\line(1,1){10}}
\put(248,22){\vector(1,1){10}}
\put(257,10){\vector(-1,1){22}}
\put(141,0){$\varepsilon_x = +1$}
\put(232,0){$\varepsilon_x = -1$}
\put(152.5,28){$x$}
\put(245,28){$x$}
\end{picture}

As above, a chord diagram of the singular knot is the circle with pre-images
of double points connected with  chords.

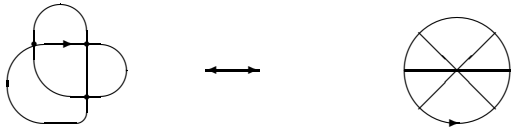
\begin{picture}(400,65)
\put(100,45){\vector(1,0){5}}
\put(100,50){\oval(20,20)[t]}
\put(110,30){\line(0,1){20}}
\put(105,30){\oval(50,30)[tl]}
\put(110,40){\oval(40,30)[bl]}
\put(110,25){\line(1,0){5}}
\put(95,20){\oval(30,10)[br]}
\put(115,35){\oval(20,20)[r]}
\put(95,30){\oval(30,30)[bl]}

\put(165,35){\vector(-1,0){10}}
\put(165,35){\vector(1,0){10}}
\put(250,35){\circle{40}}
\put(250,15){\vector(1,0){1}}
\put(230,35){\line(1,0){40}}
\put(235.5,20.5){\line(1,1){29}}
\put(264.5,20.5){\line(-1,1){29}}

\put(110,25){\circle*{2}}
\put(110,20){\line(0,1){10}}
\put(90,45){\circle*{2}}
\put(110,45){\circle*{2}}
\put(90,40){\line(0,1){10}}
\put(105,45){\line(1,0){10}}

\end{picture}

\hskip3cm {$K^\sing_3$ \qquad \qquad \qquad \qquad \qquad \qquad $D_3$}

To obtain the analogus diagram
of an ordinary knot (that is called an arrow diagram) from the chord diagram
of corresponding singular knot we must give the information on overpasses and
underpasses. Each chord is oriented from the upper branch to the lower one
and equipped with the sign (the local writh number of corresponding double
point of planar projection of the knot).

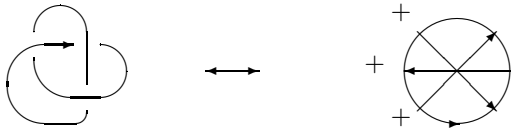
\begin{picture}(400,65)
\put(100,45){\vector(1,0){5}}
\put(100,50){\oval(20,20)[t]}
\put(110,30){\line(0,1){20}}
\put(105,30){\oval(50,30)[tl]}
\put(110,40){\oval(40,30)[bl]}
\put(110,25){\line(1,0){5}}
\put(95,20){\oval(30,10)[br]}
\put(115,35){\oval(20,20)[r]}
\put(95,30){\oval(30,30)[bl]}

\put(165,35){\vector(-1,0){10}}
\put(165,35){\vector(1,0){10}}
\put(250,35){\circle{40}}
\put(250,15){\vector(1,0){1}}
\put(270,35){\vector(-1,0){40}}
\put(235,20){\vector(1,1){30}}
\put(235,50){\vector(1,-1){30}}
\put(215,35){$+$}
\put(225,15){$+$}
\put(225,55){$+$}
\end{picture}

\hskip3cm {$K$\qquad \qquad \qquad \qquad \qquad \qquad \qquad $G$}

\section{Formulas for Vassiliev invariants of orders 2, 3 and 4.}
\subsection{Formulas of Lannes.}
We present formulas of Lannes in new form.
$$V_2(K)=1/2 \sum_{\{x,y\}\in P_2}(-1)^{\delta_x + \delta_y} W_2(\{x,y\})
\varepsilon_x \varepsilon_y [\delta_x (1- \delta_y )+\delta_y (1-\delta_x)],$$
$$V_3(K)=1/2 \sum_{\{x,y,z\}\in P_3} (-1)^{\delta_x+\delta_y+\delta_z}
W_3(\{x,y,z\}) \varepsilon_x \varepsilon_y \varepsilon_z
[\delta_y (1- \delta_x)(1- \delta_z)- \delta_x \delta_z (1- \delta_y)],$$
where the sum is taken over all unordered pairs (triplets) of double points
of planar projection, $W_2(\{x,y\})$ ($W_3(\{x,y,z\})$) is weight of chord
diagram corresponding to pair (triplet) of double points.

\subsection{Formulas of Viro-Polyak.}
Denote by $<A,G>$ algebraic number of subdiagrams of given combinatoial type
$A \quad A \subset G$ and let $<\sum_{i}n_iA_i,G> = \sum_{i}n_i<A_i,G>,
\quad n_i \in {\bf Q}$ by definition. Then

\begin{picture}(400,22)
\put(193,5){\vector(1,1){11}}
\put(203,5){\vector(-1,1){11}}
\put(198,10){\circle{15}}
\put(198,17){\circle*{2}}
\put(135,5){\it $V_2(K)=<$}
\put(210,5){$,G>,$}
\end{picture}

\begin{picture}(400,22)
\put(183,5){\vector(1,1){11}}
\put(193,5){\vector(-1,1){11}}
\put(188,17){\vector(0,-1){15}}
\put(188,10){\circle{15}}
\put(115,5){\it $V_3(K)=<[$}
\put(260,5){$],G>,$}
\put(205,5){$]+\frac{1}{2}[$}
\put(243,13){\vector(1,0){11}}
\put(253,5){\vector(-1,0){11}}
\put(248,2){\vector(0,1){15}}
\put(248,10){\circle{15}}
\end{picture}

\begin{picture}(420,20)
\put(-10,5){\it $V_4(K)=<$}
\put(145,5){$+2$}
\put(173,8){\circle{15}}
\put(170,2){\vector(0,1){12}}
\put(176,14){\vector(0,-1){12}}
\put(181,8){\vector(-1,0){15}}
\put(173,0){\vector(0,1){15}}
\put(185,5){-}
\put(203,8){\circle{15}}
\put(200,2){\vector(0,1){12}}
\put(206,2){\vector(0,1){12}}
\put(196,8){\vector(1,0){15}}
\put(203,16){\vector(0,-1){15}}

\put(48,8){\circle{15}}
\put(60,5){+2}
\put(45,14){\vector(0,-1){12}}
\put(51,2){\vector(0,1){12}}
\put(55,6){\vector(-1,0){14}}
\put(41,10){\vector(1,0){14}}

\put(88,8){\circle{15}}
\put(100,5){+6}
\put(85,14){\vector(0,-1){12}}
\put(91,2){\vector(0,1){12}}
\put(81,5){\vector(1,0){14}}
\put(95,11){\vector(-1,0){14}}

\put(128,8){\circle{15}}
\put(125,2){\vector(0,1){12}}
\put(131,14){\vector(0,-1){12}}
\put(135,5){\vector(-1,0){14}}
\put(121,11){\vector(1,0){14}}

\put(233,8){\circle{15}}
\put(233,0){\vector(-1,2){6}}
\put(233,15){\vector(-1,-2){6}}
\put(240,6){\vector(-3,-1){11}}
\put(229,15){\vector(2,-1){12}}
\put(243,5){$+$}
\put(213,5){$+$}
\put(263,8){\circle{15}}
\put(263,0){\vector(-1,2){7}}
\put(257,4){\vector(1,2){6}}
\put(269,5){\vector(-3,-1){11}}

\put(259,15){\vector(2,-1){12}}

\put(275,5){$+$}
\put(293,8){\circle{15}}
\put(300,4){\vector(-1,0){12}}
\put(286,9){\vector(1,0){11}}
\put(286,9){\vector(1,0){14}}
\put(293,0){\vector(0,1){15}}

\put(305,5){$+$}
\put(328,8){\circle{15}}
\put(328,16){\vector(0,-1){15}}
\put(328,16){\vector(0,-1){13}}
\put(321,5){\vector(2,1){14}}
\put(335,5){\vector(-2,1){14}}

\put(340,5){$+2$}
\put(368,8){\circle{15}}
\put(365,14){\vector(1,-2){6}}
\put(370,14){\vector(-1,-2){6}}
\put(360,4){\vector(2,1){14}}
\put(374,4){\vector(-2,1){14}}
\put(383,5){$+$}
\end{picture}

\begin{picture}(420,20)
\put(58,8){\circle{15}}
\put(70,5){+}
\put(59,0){\vector(-1,2){7}}
\put(59,16){\vector(-1,-2){7}}
\put(61,14){\vector(0,-1){12}}
\put(51,8){\vector(1,0){15}}

\put(88,8){\circle{15}}
\put(100,5){+2}
\put(83,14){\vector(1,-2){6}}
\put(83,2){\vector(1,2){6}}
\put(92,14){\vector(0,-1){12}}
\put(96,8){\vector(-1,0){16}}

\put(128,8){\circle{15}}
\put(129,0){\vector(-1,2){7}}
\put(123,2){\vector(1,2){6}}
\put(132,14){\vector(0,-1){12}}
\put(136,8){\vector(-1,0){16}}

\put(140,5){+2}
\put(173,8){\circle{15}}
\put(170,2){\vector(0,1){12}}
\put(176,2){\vector(0,1){12}}
\put(166,5){\vector(2,1){14}}
\put(180,5){\vector(-2,1){14}}
\put(185,5){+3}
\put(213,8){\circle{15}}
\put(210,2){\vector(0,1){12}}
\put(216,14){\vector(0,-1){12}}
\put(206,5){\vector(2,1){14}}
\put(219,4){\vector(-2,1){14}}
\put(225,5){+}
\put(256,5){$,G>$}
\put(243,8){\circle{15}}
\put(240,14){\vector(0,-1){12}}
\put(245,2){\vector(0,1){12}}
\put(235,4){\vector(2,1){14}}
\put(249,4){\vector(-2,1){14}}
\end{picture}

\subsection{New formula.}
\noindent
\begin{picture}(420,20)
\put(10,5){\it $V_4(K)=\frac{1}{2}W_4($}
\put(93,8){\circle{15}}      
\put(91,1){\line(-1,2){5}}   
\put(86,4){\line(1,2){5.5}}  
\put(95,1){\line(1,2){5}}    
\put(100,4){\line(-1,2){5.5}}  
\put(102,6){$)<$}
\put(128,8){\circle{15}}      
\put(122,12){\vector(1,-2){6}}   
\put(121,4){\vector(1,2){5.5}}    
\put(130,1){\vector(1,2){5.5}}    
\put(129,15){\vector(1,-2){6}}   
\put(140,5){$,G>+\frac{1}{4}W_4($}
\put(213,8){\circle{15}}      
\put(210,2){\line(0,1){12}}   
\put(216,2){\line(0,1){12}}   
\put(206,8){\line(1,0){15}}   
\put(213,0){\line(0,1){15}}   
\put(225,5){$)<$}
\put(253,8){\circle{15}}      
\put(250,14){\vector(0,-1){12}}   
\put(256,14){\vector(0,-1){12}}   
\put(246,8){\vector(1,0){15}}   
\put(253,0){\vector(0,1){15}}   
\put(265,5){+}
\put(283,8){\circle{15}}        
\put(280,2){\vector(0,1){12}}   
\put(286,2){\vector(0,1){12}}   
\put(276,8){\vector(1,0){15}}   
\put(283,16){\vector(0,-1){15}} 
\put(300,5){$,G>+$}
\end{picture}

\hskip1.2cm
\begin{picture}(420,20)
\put(5,5){$\frac{1}{2}W_4($}
\put(43,8){\circle{15}}     
\put(52,5){$)<$}           
\put(40,2){\line(0,1){12}}  
\put(46,2){\line(0,1){12}}  
\put(36,5){\line(1,0){14}}  
\put(36,11){\line(1,0){14}} 
\put(78,8){\circle{15}}     
\put(90,5){+}               
\put(75,14){\vector(0,-1){12}}  
\put(81,2){\vector(0,1){12}}  
\put(71,5){\vector(1,0){14}}  
\put(85,11){\vector(-1,0){14}} 
\put(108,8){\circle{15}}     
\put(120,5){$,G>+\frac{1}{4}W_4($}           
\put(105,2){\vector(0,1){12}}  
\put(111,14){\vector(0,-1){12}}  
\put(115,5){\vector(-1,0){14}}  
\put(101,11){\vector(1,0){14}} 
\put(198,8){\circle{15}}     
\put(210,5){$)<$}            
\put(198,0){\line(-1,2){6}}  
\put(192,4){\line(1,2){6}}   
\put(194,2){\line(3,1){11}}  
\put(194,15){\line(2,-1){12}}
\put(238,8){\circle{15}}         
\put(250,5){+}                   
\put(238,0){\vector(-1,2){6}}    
\put(238,15){\vector(-1,-2){6}}  
\put(245,6){\vector(-3,-1){11}}  
\put(234,15){\vector(2,-1){12}}  
\put(268,8){\circle{15}}         
\put(280,5){+}                   

\put(263,14){\vector(1,-2){7}}    
\put(262,4){\vector(1,2){6}}  
\put(274.5,6){\vector(-3,-1){11}}  
\put(264,15){\vector(2,-1){12}}  

\put(328,8){\circle{15}}         
\put(340,5){,$G>+$}             

\put(328,0){\vector(-1,2){6}}    
\put(328,15){\vector(-1,-2){6}}  
\put(324,1){\vector(3,1){11}}  
\put(336,8.5){\vector(-2,1){12.5}}  

\put(298,8){\circle{15}}         
\put(310,5){+}                   
\put(292,12){\vector(1,-2){6}}    
\put(293,3){\vector(1,2){6}}  
\put(293.5,1.5){\vector(3,1){11}}  
\put(304.5,9){\vector(-2,1){11.5}}  
\end{picture}

\hskip1.2cm
\begin{picture}(420,20)
\put(5,5){$\frac{1}{4}W_4($}
\put(43,8){\circle{15}}     
\put(52,5){$)<$}           
\put(43,0){\line(-1,2){6}}  
\put(37,4){\line(1,2){6}}   
\put(46,2){\line(0,1){12}}  
\put(36,8){\line(1,0){15}}  
\put(78,8){\circle{15}}     
\put(90,5){+}              
\put(79,0){\vector(-1,2){6.5}}  
\put(79,16){\vector(-1,-2){6.5}}   
\put(81,2){\vector(0,1){12}}  
\put(71,8){\vector(1,0){15}}  
\put(108,8){\circle{15}}     
\put(120,5){+}              
\put(109,0){\vector(-1,2){6.5}}  
\put(109,16){\vector(-1,-2){6.5}}   
\put(111,14){\vector(0,-1){12}}  
\put(101,8){\vector(1,0){15}}  
\put(138,8){\circle{15}}     
\put(150,5){+}              
\put(133,14){\vector(1,-2){6.5}}  
\put(133,2){\vector(1,2){6.5}}   
\put(142,14){\vector(0,-1){12}}  
\put(146,8){\vector(-1,0){16}}  
\put(168,8){\circle{15}}     
\put(180,5){,$G>+\frac{1}{4}W_4($} 
\put(163,14){\vector(1,-2){6.5}}  
\put(163,2){\vector(1,2){6.5}}   
\put(172,2){\vector(0,1){12}}  
\put(176,8){\vector(-1,0){16}}  
\put(253,8){\circle{15}}     
\put(250,2){\line(0,1){12}}  
\put(256,2){\line(0,1){12}}  
\put(246,5){\line(2,1){14}}  
\put(246,11){\line(2,-1){14}}
\put(265,5){$)<$}
\put(293,8){\circle{15}}     
\put(290,2){\vector(0,1){12}}  
\put(296,14){\vector(0,-1){12}}  
\put(286,5){\vector(2,1){14}}  
\put(299,4){\vector(-2,1){13.54}}
\put(305,5){+}
\put(323,8){\circle{15}}     
\put(320,14){\vector(0,-1){12}}  
\put(325,2){\vector(0,1){12}}  
\put(316,4){\vector(2,1){14}}  
\put(329,4){\vector(-2,1){13.54}}
\put(335,5){,$G>$}
\end{picture}

\noindent
{\it for any weight system $W_4$ such that $W_4(
\begin{picture}(20,15)
\put(8,3){\circle{15}}     
\put(3,-2){\line(1,1){10}}  
\put(13,-2){\line(-1,1){10}} 
\put(8,-5){\line(0,1){15}}  
\put(1,3){\line(1,0){15}}  
\end{picture}
)=0$}

\end{document}